\documentclass[leqno]{amsart}
\usepackage{fullpage}
\usepackage[T1]{fontenc}
\usepackage[utf8]{inputenc}
\usepackage[english]{babel}
\usepackage{amsmath}
\usepackage{amsfonts}
\usepackage{amssymb}
\usepackage{amsthm}
\usepackage{float}
\usepackage{mathtools}
\usepackage{enumitem}
\usepackage[colorlinks=true, allcolors=blue]{hyperref}
\hypersetup{
	colorlinks=true,
	linkcolor=blue,
	filecolor=blue,      
	urlcolor=magenta,
	citecolor=magenta,
    pdfpagemode=FullScreen,
    pdftitle={Parallel transport Kummer IHS}}

{\theoremstyle{plain}
\newtheorem{theo}{Theorem}[section]
\newtheorem{lem}[theo]{Lemma}
\newtheorem{coro}[theo]{Corollary}
\newtheorem{propo}[theo]{Proposition}
}

{%
\theoremstyle{definition}
\newtheorem{defin}[theo]{Definition}
\newtheorem{rem}[theo]{Remark}
\newtheorem*{notat}{Notation}
}

\title{Note on parallel transport of IHS manifolds of Kummer type}
\author{Simon Brandhorst and Stevell Muller}
\begin{document}
\fontsize{10}{13}
\begin{abstract}
We give a characterization of parallel transport operators for irreducible holomorphic symplectic manifolds of Kummer type.
\end{abstract}
\maketitle

\section{Introduction}
The Hodge-theoretic global Torelli theorem of Huybrechts--Markman--Verbitsky is a key tool to study irreducible holomorphic symplectic (IHS) manifolds via lattice theory. 
Let $X$ and $Y$ be two IHS manifolds which are deformation equivalent. The aforementioned theorem implies that an invertible linear map $g\colon H^2(Y,\mathbb{Z)} \to H^2(X,\mathbb{Z})$ is induced by a biholomorphic (resp. bimeromorphic) map $f\colon X \dashrightarrow Y$, i.e. $g=f^*$ if and only if $g$ preserves the weight two Hodge structures, the Kähler cones (resp. fundamental exceptional chambers), the Beauville--Bogomolov--Fujiki (BBF) quadratic forms and if $g$ is a parallel transport operator. The latter means that $g$ is induced by a parallel transport in the local system $R^2\pi_\ast\mathbb{Z}$ along a path $\gamma$, where $\pi\colon \mathcal{X}\to B$ is a smooth and proper family of IHS manifolds over an analytic base $B$ and $\gamma\colon [0,1]\to B$ is a continuous path connecting $\mathcal{X}_{\gamma(0)}\cong X$ and $\mathcal{X}_{\gamma(1)}\cong Y$.
Thus, in order to leverage such a Torelli-type theorem to compare IHS manifolds, one needs to determine when such a map $g$ is a parallel transport operator. We have a characterization for $3$ of the $4$ known deformation types of IHS manifolds \cite[Theorem~9.8]{mar11}, \cite{mr21}, \cite{ono22}. In this work we resolve the missing case of Kummer type IHS manifolds. 

With parallel transport under control, one can streamline proofs of numerical characterizations of the geometry of Kummer type IHS manifolds, e.g. when a given Kummer type manifold is an actual generalized Kummer manifold (and not just a deformation of one), or when it comes from a moduli space of (twisted) sheaves of an abelian surface (see for instance \cite[Proposition~2.16]{dmmn26}). In upcoming works, we use it to complete the classification of finite group actions on Kummer type IHS manifolds. In view of our applications, we take care to formulate the results for $\operatorname{Kum}_n$-type and $\operatorname{K3}^{[n]}$-type IHS manifolds uniformly. 

\begin{rem}
    Our result relies on the fact that the monodromy group of $\operatorname{Kum}_n$-type IHS manifolds is known (see \cite[Theorem~4.3]{mon16a} and \cite[Theorem~1.4]{mar23}), and we therefore do not prove it. It would be interesting to know whether our approach could be used to find a more natural proof for the determination of the monodromy groups of Kummer type IHS manifolds.
\end{rem}

 \section{Preliminaries}

\subsection*{Lattice notation}
In this note, a lattice is a finitely generated free $\mathbb{Z}$-module $L$ equipped with a nondegenerate $\mathbb{Z}$-valued symmetric bilinear form $b$ whose associated quadratic form $x \mapsto b(x,x)$ takes even integer values (i.e., the lattice is often said to be \emph{even}). If there is no ambiguity, we write $x^2 \coloneqq b(x,x)$ and $x.y \coloneqq b(x,y)$ for all $x,y\in L$. We always assume that the ADE root lattices associated to the corresponding Dynkin diagrams are negative definite, and $U$ denotes the hyperbolic plane lattice.\smallskip

Orthogonal direct sums of lattices are denoted using the symbol $\oplus$.\smallskip

Given a lattice $L$, we denote by $O(L)$ its group of isometries and by ${\rm{SO}}(L)$ the \emph{special orthogonal group}, consisting of isometries of positive determinant. 
We let $D_L\coloneqq L^\vee/L$ denote its \emph{discriminant group}. For every $f\in O(L)$, we then let $D_f\in \operatorname{Aut}(D_L)$ be the induced group automorphism.\smallskip

By our convention, any morphism of lattices is injective so we talk about \emph{embeddings}. Such an embedding $S\hookrightarrow L$ is called \emph{primitive} if its cokernel is torsionfree. Given two lattices $S$ and $L$, and a subgroup $G_L\leq O(L)$ (resp. $G_S\leq O(S)$) there is a left action (resp. right action) of $G_L$ (resp. $G_S$) on the set of primitive embeddings $S\hookrightarrow L$ given by postcomposition (resp. precomposition). We can therefore talk about orbits of primitive embeddings, or double cosets of such (if both a left and a right action are prescribed).\smallskip

Given a lattice $L$ and an anisotropic vector $v\in L\otimes \mathbb{R}$, we let $\alpha_v\in O(L\otimes\mathbb{R})$ be the associated reflection and we define
\[\vartheta_L(\alpha_v) := - v^2(\mathbb{R}^\times)^2 \in \mathbb{R}^\times/(\mathbb{R}^\times)^2.\]
This defines a uniquely determined group homomorphism
\[\vartheta_L\colon O(L)\to \{\pm1\}\]
which we refer to as \emph{spinor norm morphism} \cite[54:6]{om73}. For any subgroup $G\leq O(L)$, we denote $G^+ := G\cap \ker(\vartheta_L)$.

\subsection*{Mukai lattice of K-trivial surfaces}
Let $S$ be a K-trivial surface, i.e. $S$ is a compact K\"ahler K3 surface or a complex 2-torus. For any given subring $R$ of $\mathbb{C}$, we denote 
\[H^{2\ast}(S, R) := H^0(S, R)\oplus H^2(S, R)\oplus H^4(S, R).\]
We equip $H^{2\ast}(S,\mathbb{Z})$ with the standard \emph{Mukai pairing}
\begin{equation}\label{mukai pairing}((v_0, v_1, v_2),\ (v_0', v_1', v_2'))\mapsto v_1\smile v_1' - v_0v_2' - v_2v_0'\end{equation}
where $\smile$ is the cup-product on $H^2(S, \mathbb{Z})$ (seen as a $\mathbb{Z}$-valued bilinear form via $H^4(S, \mathbb{Z})\cong\mathbb{Z})$. This gives $H^{2\ast}(S, \mathbb{Z})$ the structure of a unimodular lattice of signature $(4, b_2(S)-2)$. We refer to the lattice $H^{2\ast}(S, \mathbb{Z})$ as \emph{Mukai lattice} of $S$. Note that we have the following identifications:
\[H^{2\ast}(S, \mathbb{Z})\cong \begin{cases}
    U^{\oplus4}\oplus E_8^{\oplus2}&\text{if $S$ is a K3 surface}\\
    U^{\oplus4}&\text{if $S$ is a complex $2$-torus}
\end{cases}.\]
It is also equipped with a pure Hodge structure of weight two which is polarized by the previous bilinear form, and whose $(2,0)$-part is exactly $H^{2,0}(S)$.\smallskip

Set $\varepsilon = 0$ if $S$ is a K3 surface and $\varepsilon=1$ for $S$ a complex $2$-torus. For a choice of a primitive, positive Mukai vector $v\in H^{2\ast}(S, \mathbb{Z})$ of square $2n-2+4\varepsilon\geq 2+4\varepsilon$ and of a $v$-generic stability condition $\sigma$ on $S$, the Albanese fiber $S_n^\varepsilon$ of the moduli space of $\sigma$-stable sheaves on $S$ with Mukai vector $v$ is an IHS manifold. There is moreover a natural identification of integral Hodge structures
\begin{equation}\label{eq:mukai embed}
    H^2(S_n^\varepsilon,\mathbb{Z})\cong v^\perp\subset H^{2\ast}(S, \mathbb{Z}).
\end{equation}
The associated primitive embedding $\theta_v\colon H^2(S_n^\varepsilon, \mathbb{Z})\hookrightarrow H^{2\ast}(S, \mathbb{Z})$ is referred to as \emph{Mukai's homomorphism}.

\subsection*{Abstract Mukai lattice and numerical monodromy}
For what follows, we set $\varepsilon=0$ for the $\text{K3}^{[n]}$-cases, and $\varepsilon=1$ for the $\text{Kum}_n$-cases. Some definition and notation depend on the choice of such $\varepsilon$; we try as much as possible to avoid any ambiguity. For instance, given a lattice $L$, we denote
\[S(L) \coloneqq \left\{\begin{array}{ccc}O(L)&\text{ if }&\varepsilon=0\\
{\rm{SO}}(L)&\text{ if }&\varepsilon=1\end{array}\right..\]
Moreover, given any $n\geq 2$, we let
\[\Lambda_n^{\varepsilon}\coloneqq U^{\oplus3}\oplus E_8^{\oplus (2-2\varepsilon)}\oplus A_1(n-1+2\varepsilon)\]
be the usual reference lattice associated to the corresponding deformation type, and we let
\[M^{\varepsilon} \coloneqq U^{\oplus4}\oplus E_8^{\oplus (2-2\varepsilon)}\]
be the associated abstract Mukai lattice.

\begin{rem}
    Here, for a given lattice $L$, we let $L^{\oplus0} = \{0\}$ be the trivial lattice.
\end{rem}

As we have already seen in (\ref{eq:mukai embed}), for any $n\geq 2$, the lattice $\Lambda_n^\varepsilon$ embeds primitively into $M^\varepsilon$. Given any primitive embedding $\kappa\colon \Lambda_n^{\varepsilon}\hookrightarrow M^{\varepsilon}$, we denote by $V_\kappa = \mathbb{Z}v_\kappa$ where $(v_\kappa)^2 = 2(n-1+2\varepsilon)$, the associated orthogonal complement. We let moreover $S_\kappa\coloneqq S(M^\varepsilon, V_\kappa)$ be the setwise stabilizer of $V_\kappa$ in $S(M^\varepsilon)$. We define $\vartheta_\kappa\colon S_\kappa\to \{\pm1\}$ to be the restriction of the spinor norm morphism $\vartheta_{M^\varepsilon}$ to $S_\kappa$ and $\chi_\kappa\colon S_\kappa\to \{\pm1\}$ to be the character such that $V_\kappa\leq \ker(g-\chi_\kappa(g)\cdot\text{id})$ for all $g\in S_\kappa$.
Since isometries in $S_\kappa$ preserve $V_\kappa$, they act on $\Lambda^{\varepsilon}_n$. Thus, there exists a well-defined restriction map $p_\kappa\colon S_\kappa\to O(\Lambda^{\varepsilon}_n)$. 

\begin{lem}\label{lem:abstract monodromy}
    The image of $p_\kappa$ does not depend on the primitive embedding $\kappa\colon \Lambda_n^\varepsilon \hookrightarrow M^\varepsilon$.
\end{lem}

\begin{proof}
        Any isometry in $S_\kappa$ acts by $\pm \text{id}_{V_\kappa}$ on $V_\kappa$. According to Nikulin's condition on equivariant primitive extensions (see \cite[Corollary~1.5.2]{nik80}), since $M^\varepsilon$ is unimodular, the image of $p_\kappa$ consists of all isometries $f\in O(\Lambda^{\varepsilon}_n)$ acting as $\pm 1$ on $D_{\Lambda^{\varepsilon}_n}$, with the extra condition that $D_f = \text{det}(f)\text{id}_{D_{\Lambda^{\varepsilon}_n}}$ when $\varepsilon=1$. The subgroup of $O(\Lambda^{\varepsilon}_n)$ consisting of such isometries only depends on $\Lambda^{\varepsilon}_n$.
\end{proof}

\begin{defin}
    We define the \emph{numerical monodromy group} of $\Lambda^{\varepsilon}_n$ to be 
    \[\operatorname{Mon}^2(\Lambda^{\varepsilon}_n) \coloneqq \text{im}(p_\kappa)^+ = \operatorname{im}(p_{\kappa})\cap \ker(\vartheta_{\Lambda_n^\varepsilon}),\]
    for any primitive embedding $\kappa\colon \Lambda_n^\varepsilon \hookrightarrow M^\varepsilon$.
\end{defin}

\begin{rem}\label{rem monodromy}
The subgroup $\operatorname{Mon}^2(\Lambda^{\varepsilon}_n)\leq O(\Lambda_n^\varepsilon)$ is normal of finite index. If $X$ is an IHS manifold of $\operatorname{K3}^{[n]}$-type or of $\operatorname{Kum}_n$-type, and if $\eta\colon H^2(X,\mathbb{Z})\to \Lambda_n^\varepsilon$ is any isometry, then by \cite[Lemma~4.2]{mar10}, \cite[Theorem~4.3]{mon16a} and \cite[Theorem~1.4]{mar23} the following holds:
\[\operatorname{Mon}^2(X) = \eta^{-1}\operatorname{Mon}^2(\Lambda_n^\varepsilon)\eta.\]
\end{rem}

The corestriction of $p_{\kappa}$ to $\operatorname{Mon}^2(\Lambda^{\varepsilon}_n)$ admits a
section 
\begin{equation}\label{gamma map} \gamma_\kappa\colon \textnormal{Mon}^2(\Lambda^{\varepsilon}_n)\to S_\kappa, \; h\mapsto \widehat{\chi}_\kappa(h)\oplus h\end{equation}
where we define $\widehat{\chi}_\kappa(h) = \operatorname{id}_{V_\kappa}$ if $D_h=\text{id}_{D_{\Lambda^{\varepsilon}_n}}$, and $-\operatorname{id}_{V_\kappa}$ otherwise.
According to \cite[Lemma 4.6]{bc22} the image of $\gamma_\kappa$ is given by the kernel $\Gamma_\kappa$ of
\[ \vartheta_\kappa\cdot \chi_\kappa\colon S_\kappa\to \{\pm 1\}.\]

\section{Parallel transport}
From now on, and for the rest of the paper, we fix an orientation on $M^1 = U^{\oplus4}$, i.e an ${\rm{SO}}(M^1)$-orbit of a basis. In this section we  strengthen \cite[Theorem~A]{dmmn26} by removing the conditions on \emph{compatible embeddings}, as defined in \cite[Definition~2.11]{dmmn26}.
In order to do so, we start by showing that there is a natural way of orienting the Mukai lattice of an abelian surface.

\begin{propo}\label{propo:natural orientation}
    Let $A$ be an abelian surface. Its Mukai lattice $H^{2*}(A,\mathbb{Z})$, together with the usual Mukai pairing (\ref{mukai pairing}), carries a natural orientation.
    Here by ``natural'' we mean that it only depends on the topological structure of $A$.
\end{propo}

\begin{proof}
Pick any nonzero class $x \in H^1(A, \mathbb{R})$. Then $I(x) \coloneqq \mathbb{R} [A] \oplus \left(\mathbb{R} x \wedge H^1(A,\mathbb{R})\right)$ is a maximal totally isotropic subspace of $H^{2*}(A, \mathbb{R})$. 
It defines an orientation as follows: let $(b_1,\dots, b_4)$ be any basis of $I(x)$.
Choose $b_5,\dots, b_8 \in H^{2*}(A,\mathbb{R})$ satisfying $b_i.b_{j+4}=\delta_{ij}$ for all $i,j\in\{1,2,3,4\}$ (where $\delta$ denotes the Kronecker delta). This choice is unique modulo $I(x)$ and therefore the orientation defined by the ordered basis $(b_1,\dots, b_8)$ does not depend on the choice of $b_5,\dots, b_8$. 
Changing $b_1,\dots, b_4$ by an element $f \in \mathrm{GL}(I(x))$ changes $b_5,\dots, b_8$ by $f^{-T}$ modulo $I(x)$. 
The total change of basis has determinant one. 
A priori the orientation depends on the choice of $x$. But $H^1(A,\mathbb{R}) \setminus \{0\}$ is connected, and so any two choices can be deformed into one another. Thus they define the same orientation.
\end{proof}

We can now generalize Markman's construction \cite[Corollary 9.5]{mar11} of a natural orbit of primitive embeddings into the Mukai lattice, from the $\operatorname{K3}^{[n]}$-cases to the $\operatorname{Kum}_n$-cases. Such a class of primitive embeddings was already constructed by Wieneck \cite[Theorem~4.9]{wie18} in a coarser sense.

\begin{propo}\label{prop:natural orbit of embeddings}
    Let $n\geq 2$ and let $X$ be IHS manifold which is of deformation type
    \begin{enumerate}
        \item[(a)] $\textnormal{K3}^{[n]}$-type if $\varepsilon=0$, or
        \item[(b)] $\textnormal{Kum}_{n}$-type if $\varepsilon=1$.
    \end{enumerate}
    Then there exists a natural $S(M^\varepsilon)$-orbit of primitive embeddings $\kappa_X\colon H^2(X, \mathbb{Z})\hookrightarrow M^\varepsilon$. 
\end{propo}
\begin{proof}
For the $K3^{[n]}$-type we refer to \cite[Corollary 9.5]{mar11}.

By definition we can deform $X$ to a generalized Kummer manifold $Y= \textnormal{K}_n(A)$ of an abelian surface $A$. 
Pick any orientation-preserving isometry $\eta \colon H^{2*}(A,\mathbb{Z}) \to M^1$, for the fixed orientation on $M^1$ and the natural orientation on $H^{2*}(A,\mathbb{Z})$ discussed above.
By composing some parallel transport $P_\alpha$ along some path $\alpha$ in the deformation and Mukai's homomorphism $\theta_v$ we obtain a primitive embedding:
\[\kappa_{\alpha,\eta} \colon H^2(X,\mathbb{Z})\xrightarrow[]{P_\alpha} H^2(Y, \mathbb{Z})\xrightarrow[\cong]{\theta_v} v^\perp \hookrightarrow H^{2*}(A,\mathbb{Z}) \xrightarrow[]{\eta} M.\]
Changing the path $\alpha$ changes the embedding $\kappa_{\alpha,\eta}$ by a monodromy $m \in \operatorname{Mon}^2(X)$, i.e., to
\[\kappa_{\alpha, \eta} \circ m = \gamma_{\kappa_{\alpha,\eta}}(m) \circ \kappa_{\alpha,\eta},\]
where $\gamma_{\kappa_{\alpha,\eta}}$ is defined as in (\ref{gamma map}). Since $\gamma_{\kappa_{\alpha,\eta}}(m) \in {\rm{SO}}(M^1)$, this does not change the orbit. Since $\eta$ is orientation-preserving, any two different choices for $\eta$ differ by an element of ${\rm{SO}}(M^1)$ as well. Finally, if we choose a different abelian surface $A'$, then we can deform $A$ to $A'$ and this deformation is compatible with Mukai's homomorphism $\theta_v$ by \cite[Proposition 5.1]{yos01}. Further, parallel transport along this deformation preserves the natural orientation on $H^2(A,\mathbb{Z})$, by the proof of Proposition~\ref{propo:natural orientation}.\qedhere
\end{proof}

Following \cite[Section 4]{mar11} we recall the notion of an orientation of a lattice $\Lambda$ of signature $(3,\ast)$. Let $C_\Lambda := \{x \in \Lambda \otimes \mathbb{R} \mid x^2 >0 \}$. The choice of a generator $\alpha$ of $H^2(C_\Lambda,\mathbb{Z}) \cong \mathbb{Z}$ is called an \emph{orientation} of $\Lambda$. Such an orientation is determined by any oriented, positive definite, 3-dimensional subspace $W \subseteq \Lambda_\mathbb{R}$. The group of isometries preserving such an orientation coincides with $O^+(L)$ \cite[Lemma~2.2]{bmm26}.

Given an IHS manifold $X$, the lattice $H^2(X,\mathbb{Z})$ comes equipped with a natural orientation, in the above sense, defined by the oriented, positive definite, real quadratic space
$\operatorname{span}_{\mathbb{R}}(\text{Re}~\sigma_X, \text{Im}~\sigma_X, h)$ where $h$ is any K\"ahler class on $X$. 

\begin{notat}
    In order to avoid confusions with the natural orientation of the Mukai lattice of an abelian surface, we refer to the natural orientation on the BBF-lattice of an IHS manifold as \emph{Markman's orientation}.
\end{notat}

\begin{propo}\label{prop:extension of parallel transport}
    Let $n\geq 2$ and let $X, Y$ be IHS manifolds which are both of 
    \begin{enumerate}
        \item[(a)] $\textnormal{K3}^{[n]}$-type if $\varepsilon=0$, or
        \item[(b)] $\textnormal{Kum}_{n}$-type if $\varepsilon=1$.
    \end{enumerate}
    Then an isometry $g\colon H^2(X,\mathbb{Z})\to H^2(Y, \mathbb{Z})$ is a parallel transport operator if and only if $g$ respects Markman's orientations and the natural primitive embeddings $S(M^\varepsilon)\kappa_X = S(M^\varepsilon)(\kappa_Y \circ g)$, as defined in Proposition~\ref{prop:natural orbit of embeddings}.
\end{propo}

\begin{proof}
The $\textnormal{K3}^{[n]}$-case is \cite[Theorem 9.8]{mar11}. We prove the $\textnormal{Kum}_n$-case. Already note that any parallel transport operator $H^2(X,\mathbb{Z})\to H^2(Y, \mathbb{Z})$ respects Markman's orientations. 
Let now $g=P_\gamma$ be a parallel transport operator. Then $\kappa_X$ and $\kappa_Y\circ g$ are both defined by some parallel transport from $X$ to some (fixed) generalized Kummer manifold $K$, followed by the associated Mukai's homomorphism and a (fixed) orientation-preserving isometry.
We see that $h\coloneqq\kappa_X^{-1} \circ \kappa_Y \circ g$ is a monodromy operator and 
$\gamma_{\kappa_X}(h) \circ \kappa_X = \kappa_X \circ h =  \kappa_Y \circ g$. 

Suppose now that $g$ respects Markman's orientations and that there exists $h \in {\rm{SO}}(M^1)$ with $\kappa_Y \circ g = h \circ \kappa_X$. Since $X$ and $Y$ are deformation equivalent, there exists a parallel transport operator $f \colon H^2(Y,\mathbb{Z})\to H^2(X,\mathbb{Z})$.
It satisfies $h' \circ \kappa_Y = \kappa_X \circ f$ for some $h' \in {\rm{SO}}(M^1)$.
Then the equality $(h'h) \circ \kappa_X = \kappa_X \circ f \circ g$ shows that $h'' \coloneqq h'h \in {\rm{SO}}(M^1)$ preserves $V_{\kappa_X}$ and $f\circ g = p_{\kappa_X}(h'')$. Moreover, since $f\circ g$ preserves the Markman's orientation on $H^2(X,\mathbb{Z})$ by our assumption on $g$, the composition $f\circ g$ is a monodromy operator (see Remark~\ref{rem monodromy}). Now we can write 
$g = f^{-1} \circ (f\circ g)$ as a composition of parallel transport operators, so it is itself a parallel transport operator. 
\end{proof}

\begin{rem}\label{rem:some remark}
  The monodromy group $\operatorname{Mon}^2(X)$ of $X$ coincides with the stabilizer of $S(M^\varepsilon)\kappa_X$ for the right action of $O^+(H^2(X, \mathbb{Z}))$ via $S(M^\varepsilon)\kappa_X\cdot g \coloneqq S(M^\varepsilon)(\kappa_X\circ g)$.
\end{rem}

We conclude with the following corollary, whose proof is essentially contained in \cite[Corollary~9.9]{mar11} and \cite[Theorem~2.14]{dmmn26}.

\begin{coro}
     Let $n\geq 2$ and let $X, Y$ be IHS manifolds which are both of 
    \begin{enumerate}
        \item[(a)] $\textnormal{K3}^{[n]}$-type if $\varepsilon=0$, or
        \item[(b)] $\textnormal{Kum}_{n}$-type if $\varepsilon=1$.
    \end{enumerate}
    Then, the manifolds $X$ and $Y$ are bimeromorphic if and only if there exists a Hodge isometry $g\colon H^2(X,\mathbb{Z})\to H^2(Y, \mathbb{Z})$ such that $S(M^\varepsilon)\kappa_X = S(M^\varepsilon)(\kappa_Y\circ g)$.
\end{coro}

\section*{Acknowledgment}
We would like to thank Dominique Mattei for the various helpful discussions on the content of the note.

\section*{Disclosures}

\subsection*{Funding}
The first named author is funded by the Deutsche Forschungsgemeinschaft (DFG, German Research Foundation) – Project-ID 286237555 – TRR 195. Gefördert durch die Deutsche Forschungsgemeinschaft (DFG) – Projektnummer 286237555 – TRR 195.

\subsection*{Tool and computational resource}
Generative AI was used for final proofreading of the manuscript. It suggested a simplification of the proof of Proposition~\ref{propo:natural orientation} which we adopted. 

\subsection*{Conflict of interests}
The authors declare no conflict of interests.

\end{document}